\documentclass{article}
\usepackage{latexsym}

\begin{document}
This is a preprint of a paper whose final and definite form will be
published in Thermal Science. Paper Submitted 28/ Dec /2016; Revised
20/Jan/2016; Accepted for publication 21/Jan/2016.

\begin{center}
\textbf{Some new applications for heat and fluid flows via fractional
derivatives without singular kernel}
\end{center}

\begin{center}
\textbf{Xiao-Jun Yang}
\end{center}

\begin{center}
\textbf{State Key Laboratory for Geomechanics and Deep Underground Engineering, School of Mechanics and Civil Engineering, China University of Mining and Technology, Xuzhou, 221116, China}
\end{center}
\begin{center}
\textbf{Zhi-Zhen Zhang}
\end{center}

\begin{center}
\textbf{State Key Laboratory for Geomechanics and Deep Underground Engineering, School of Mechanics and Civil Engineering, China University of Mining and Technology, Xuzhou, 221116, China}
\end{center}
\begin{center}
\textbf{H. M. Srivastava}
\end{center}

\begin{center}
\textbf{Department of Mathematics and Statistics, University of Victoria, Victoria, British Columbia V8W 3R4, Canada}
\end{center}
\textbf{Abstract: }

This paper addresses the mathematical models for the heat-conduction
equations and the Navier-Stokes equations via fractional derivatives without
singular kernel.


\textbf{Introduction }

Fractional derivatives of variable order [1-4] has used to set up the
mathematical models for engineering practice, especially in the fields of
the heat [5-6]and fluid flows [7-8].

Recently, Caputo and Fabrizio reported the fractional derivative operator
without singular kernel, which was given as (see[9-13]):
\begin{equation}
\label{eq1}
D_x^{\left( \beta \right)} \Xi \left( x \right)=\frac{\left( {2-\beta }
\right)\aleph \left( \beta \right)}{2\left( {1-\beta }
\right)}\int\limits_0^x {\exp \left( {-\frac{\beta }{1-\beta }\left(
{x-\lambda } \right)} \right)\Xi ^{\left( 1 \right)}\left( \lambda
\right)d\lambda } ,
\end{equation}
where $\aleph \left( \beta \right)$ is a normalization constant depending on
$\beta \mbox{ }\left( {0<\beta <1} \right)$ such that $\aleph \left( 0
\right)=\aleph \left( 1 \right)=1$.

More recently, Yang, Srivastava and Machado reported a new fractional
derivative without singular kernel (see [12,13])
\begin{equation}
\label{eq2}
D_{a^+}^{\left( \beta \right)} \Omega \left( x \right)=\frac{\aleph \left(
\beta \right)}{1-\beta }\frac{d}{dx}\int\limits_a^x {\exp \left(
{-\frac{\beta }{1-\beta }\left( {x-\lambda } \right)} \right)\Omega \left(
\lambda \right)d\lambda } ,
\end{equation}
where $a\le x$, $\beta \mbox{ }\left( {0<\beta <1} \right)$ is a real
number, and $\aleph \left( \beta \right)$ is a normalization function
depending on $\beta $ such that $\aleph \left( 0 \right)=\aleph \left( 1
\right)=1$.

In this article, our aim is to set up the heat-conduction equation and the
Navier-Stokes equation via fractional derivatives without singular kernel.
The structure of the paper is presented as follows. In Section 2, the basic
theory of the fractional derivatives without singular kernel are presented.
In Section 3, the
heat-conduction equations with the fractional derivatives without singular
kernel are proposed. In Section 4, the Navier-Stokes equations with the
fractional derivatives without singular kernel are discussed. Finally, the
conclusions are outlined in Section 5.

\textbf{Mathematical tools}

The fractional gradient operator via the Caputo and Fabrizio fractional
derivative without singular kernel is given by[9]
\begin{equation}
\label{eq3}
{ }_{CF}\nabla ^{\left( \beta \right)}\wp \left( x \right)=:\frac{\beta
}{\left( {1-\beta } \right)\sqrt {\pi ^\beta } }\int\limits_\Omega {\nabla
\wp \left( x \right)\exp \left( {-\left( {\frac{\beta }{1-\beta }\left(
{x-y} \right)} \right)^2} \right)} \mbox{d}y,
\end{equation}
where $x,y\in {\Omega }$.

The fractional tensor via the Caputo and Fabrizio fractional derivative
without singular kernel is given by [9]
\begin{equation}
\label{eq4}
{ }_{CF}\nabla ^{\left( \beta \right)}A\left( x \right)=:\frac{\beta
}{\left( {1-\beta } \right)\sqrt {\pi ^\beta } }\int\limits_\Omega {\nabla
A\left( x \right)\exp \left( {-\left( {\frac{\beta }{1-\beta }\left( {x-y}
\right)} \right)^2} \right)} \mbox{d}y,
\end{equation}
where $x,y\in {\Omega }$.

The fractional Laplacian operator via the Caputo and Fabrizio fractional
derivative without singular kernel is given by [9]
\begin{equation}
\label{eq5}
{ }_{CF}\nabla ^{\left( {2\beta } \right)}\wp \left( x \right)=:\frac{\beta
}{\left( {1-\beta } \right)\sqrt {\pi ^\beta } }\int\limits_\Omega {\nabla
^2\wp \left( x \right)\exp \left( {-\left( {\frac{\beta }{1-\beta }\left(
{x-y} \right)} \right)^2} \right)} \mbox{d}y,
\end{equation}
where $x,y\in {\Omega }$.

The fractional curl via the Caputo and Fabrizio fractional derivative
without singular kernel is given by
\begin{equation}
\label{eq6}
{ }_{CF}\nabla ^{\left( \beta \right)}\times A\left( x \right)=:\frac{\beta
}{\left( {1-\beta } \right)\sqrt {\pi ^\beta } }\int\limits_\Omega {\nabla
\times A\left( x \right)\exp \left( {-\left( {\frac{\beta }{1-\beta }\left(
{x-y} \right)} \right)^2} \right)} \mbox{d}y,
\end{equation}
where $x,y\in {\Omega }$.

The fractional gradient operator via the new fractional derivative without
singular kernel is defined as:
\begin{equation}
\label{eq7}
\nabla ^{\left( \beta \right)}\wp \left( x \right)=:\nabla \left[
{\frac{\beta \Im \left( {\beta ,\mbox{n}} \right)}{\left( {1-\beta }
\right)\sqrt {\pi ^\beta } }\int\limits_\Omega {\wp \left( x \right)\exp
\left( {-\left( {\frac{\beta }{1-\beta }\left( {x-y} \right)} \right)^2}
\right)} \mbox{d}y} \right],
\end{equation}
where $x,y\in {\Omega }\in {\rm R}^n$.

We have
\begin{equation}
\label{eq8}
\begin{array}{l}
 \mathop {\lim }\limits_{\beta \to 1} \left[ {\int\limits_\Omega {\wp \left(
x \right)\exp \left( {-\left( {\frac{\beta }{1-\beta }\left( {x-y} \right)}
\right)^2} \right)} \mbox{d}y/\left( {\pi \left( {1-\beta } \right)/\beta }
\right)^{\frac{n}{2}}} \right] \\
 =\mathop {\lim }\limits_{\beta \to 1} \int\limits_\Omega {\wp \left( x
\right)\delta \left( {x-y} \right)} \mbox{d}y \\
 =\wp \left( x \right) \\
 \end{array}
\end{equation}
such that
\begin{equation}
\label{eq9}
\mathop {\lim }\limits_{\beta \to 1} \frac{\beta ^{\frac{n+2}{2}}\Im \left(
{\beta ,\mbox{n}} \right)}{\left( {1-\beta } \right)^{\frac{n+2}{2}}\pi
^{\frac{n+\beta }{2}}}=1,
\end{equation}
where $x,y\in {\Omega }\in {\rm R}^n$ and $\Im \left( {\beta ,\mbox{n}}
\right)$ is a constant.

From Eq.(\ref{eq9}), we have the following property:
\begin{equation}
\label{eq10}
\mathop {\lim }\limits_{\beta \to 1} \nabla ^{\left( \beta \right)}\wp
\left( x \right)=\nabla \wp \left( x \right).
\end{equation}
The fractional tensor via the new fractional derivative without singular
kernel is defined as:
\begin{equation}
\label{eq11}
\nabla ^{\left( \beta \right)}A\left( x \right)=:\nabla \left[ {\frac{\beta
\Im \left( {\beta ,\mbox{n}} \right)}{\left( {1-\beta } \right)\sqrt {\pi
^\beta } }\int\limits_\Omega {\wp \left( x \right)\exp \left( {-\left(
{\frac{\beta }{1-\beta }\left( {x-y} \right)} \right)^2} \right)} \mbox{d}y}
\right],
\end{equation}
where $x,y\in {\Omega }\in {\rm R}^n$ and $\Im \left( {\beta, \mbox{n}}
\right)$ is a constant.

The fractional Laplacian operator via the new fractional derivative without
singular kernel is defined as:
\begin{equation}
\label{eq12}
\nabla ^{\left( {2\beta } \right)}\wp \left( x \right)=\nabla ^2\left[
{\frac{\beta \Im \left( {\beta ,\mbox{n}} \right)}{\left( {1-\beta }
\right)\sqrt {\pi ^\beta } }\int\limits_\Omega {\wp \left( x \right)\exp
\left( {-\left( {\frac{\beta }{1-\beta }\left( {x-y} \right)} \right)^2}
\right)} \mbox{d}y} \right],
\end{equation}
where $x,y\in {\Omega }\in {\rm R}^n$ and $\Im \left( {\beta ,\mbox{n}}
\right)$ is a constant.

The fractional curl via the new fractional derivative without singular
kernel is given by
\begin{equation}
\label{eq13}
\nabla ^{\left( \beta \right)}\times A\left( x \right)=:\nabla \times \left[
{\frac{\beta \Im \left( {\beta ,\mbox{n}} \right)}{\left( {1-\beta }
\right)\sqrt {\pi ^\beta } }\int\limits_\Omega {\wp \left( x \right)\exp
\left( {-\left( {\frac{\beta }{1-\beta }\left( {x-y} \right)} \right)^2}
\right)} \mbox{d}y} \right],
\end{equation}
where $x,y\in {\Omega }\in {\rm R}^n$ and $\Im \left( {\beta ,\mbox{n}}
\right)$ is a constant.

In Eq.(\ref{eq9}), for $n=1$, $n=2$ and $n=3$ we have
\begin{equation}
\label{eq14}
\Im \left( {\beta ,\mbox{1}} \right)=\frac{\left( {1-\beta }
\right)^{\frac{3}{2}}\pi ^{\frac{1+\beta }{2}}}{\beta ^{\frac{3}{2}}},
\Im \left( {\beta ,\mbox{2}} \right)=\frac{\left( {1-\beta } \right)^2\pi
^{\frac{2+\beta }{2}}}{\beta ^2},
\quad
\Im \left( {\beta ,\mbox{3}} \right)=\frac{\left( {1-\beta }
\right)^{\frac{5}{2}}\pi ^{\frac{3+\beta }{2}}}{\beta ^{\frac{5}{2}}}.
\end{equation}
According to the expressions (\ref{eq12}) and (\ref{eq13}) , we directly have the following
properties:
\begin{equation}
\label{eq15}
\nabla ^{\left( {2\beta } \right)}\wp \left( x \right)=\nabla \cdot \nabla
^{\left( \beta \right)}\wp \left( x \right),
\quad
\mathop {\lim }\limits_{\beta \to 1} \nabla ^{\left( \beta \right)}\times
A\left( x \right)=\nabla \times A\left( x \right).
\end{equation}
In order to discuss the problems, we replace the operators ${ }_{CF}\nabla $
and $\nabla $ by ${ }^\ast \nabla $ in this article.

\textbf{The heat-conduction problem via fractional derivatives without
singular kernel}

Following the idea [14-16], the Fourier law of the heat conduction via
fractional derivatives without singular kernel is expressed by:
\begin{equation}
\label{eq16}
\left( {x,y,z,\tau } \right)=-\kappa { }^\ast \nabla ^{\left( \beta
\right)}T\left( {x,y,z,\tau } \right),
\end{equation}
where $\kappa _{ }$ denotes the thermal conductivity of the material and ${
}^\ast \nabla ^{\left( \beta \right)}$ represents the fractional gradient
operator via the fractional derivatives without singular kernel. The Fourier
law of the heat conduction via fractional derivatives without singular
kernel in one-dimensional space was discussed in [12].

The fractional heat-conduction equation with heat generation via fractional
derivatives without singular kernel is written in the form:
\begin{equation}
\label{eq17}
\kappa { }^\ast \nabla ^{\left( {2\beta } \right)}T\left( {x,y,z,\tau }
\right)-\rho c\frac{\partial T\left( {x,y,z,\tau } \right)}{\partial
t}+g\left( {x,y,z,\tau } \right)=0,
\end{equation}
where $\rho _{ }$ and $c_{ }$ are the density and the specific heat of the
material, respectively.

The fractional heat-conduction equations within fractional derivatives
without singular kernel in the two-dimensional case read:
\begin{equation}
\label{eq18}
\kappa { }^\ast \nabla ^{\left( {2\beta } \right)}T\left( {x,y,\tau }
\right)-\rho c\frac{\partial T\left( {x,y,\tau } \right)}{\partial
t}+g\left( {x,y,\tau } \right)=0,
\end{equation}
where ${ }^\ast \nabla ^{\left( {2\beta } \right)}$ is the corresponding
fractional gradient operators via the fractional derivatives without
singular kernel, $\kappa _{ }$ denotes the thermal conductivity of the
material and both $_{ }\rho _{ }$ and $c_{ }$ are the density and the
specific heat of the material, respectively.

\textbf{The Navier-Stokes equations via fractional derivatives without
singular kernel}

We now structure the fractional velocity gradient tensor in the form:
\begin{equation}
\label{eq19}
{ }^\ast \nabla ^{\left( \beta \right)}\cdot {\rm {\bf \upsilon
}}=\frac{1}{2}\left( {{\rm {\bf \Theta }}+{\rm {\bf \Theta }}^T}
\right)+\frac{1}{2}\left( {{\rm {\bf \Theta }}-{\rm {\bf \Theta }}^T}
\right)={\rm {\bf \Lambda }}+\frac{1}{2}\left( {{\rm {\bf \Theta }}-{\rm
{\bf \Theta }}^T} \right),
\end{equation}
which leads to the fractional strain rate tensor can be written as
\begin{equation}
\label{eq20}
{\rm {\bf \Lambda }}=\frac{\mbox{1}}{\mbox{2}}\left( {{ }^\ast \nabla
^{\left( \beta \right)}\cdot {\rm {\bf \upsilon }}\mbox{+}{\rm {\bf \upsilon
}}\cdot { }^\ast \nabla ^{\left( \beta \right)}} \right),
\end{equation}
where ${\rm {\bf \upsilon }}$ is the fluid velocity, ${\rm {\bf \Theta }}={
}^\ast \nabla ^{\left( \beta \right)}\cdot {\rm {\bf \upsilon }}$ and ${\rm
{\bf \Theta }}^T={\rm {\bf \upsilon }}\cdot { }^\ast \nabla ^{\left( \beta
\right)}$

Following (\ref{eq19}) and (\ref{eq20}), we can structure the linear relation of the type of
fractional Cauchy stress
\begin{equation}
\label{eq21}
{\rm {\bf J}}=-p{\rm {\bf I}}+2\mu {\rm {\bf \Lambda }}+\lambda \left( {{
}^\ast \nabla ^{\left( \beta \right)}\cdot {\rm {\bf \upsilon }}}
\right){\rm {\bf I}},
\end{equation}
where $p$ is the thermodynamic pressure, ${\rm {\bf \Lambda }}$ is the
strain rate tensor and ${\rm {\bf I}}$ is unit vector in the field, and
$\lambda $ and $\mu $ are the bulk and shear moduli of viscosity,
respectively.

Following the idea in [16], we write the continuity equation of the
fractional flow in the form:
\begin{equation}
\label{eq22}
\frac{\partial \rho }{\partial t}+{\rm {\bf \upsilon }}\cdot \left( {{
}^\ast \nabla ^{\left( \beta \right)}\rho } \right)=0,
\end{equation}
where $\rho $ represents the fluid density.

Similarly, we give the Cauchy's equation of motion of the fractional flows
\begin{equation}
\label{eq23}
\rho \frac{\partial {\rm {\bf \upsilon }}}{\partial t}={ }^\ast \nabla
^{\left( \beta \right)}\cdot {\rm {\bf J}}+\rho {\rm {\bf b}}-\rho \left(
{{\rm {\bf \upsilon }}\cdot { }^\ast \nabla ^{\left( \beta \right)}}
\right){\rm {\bf \upsilon }},
\end{equation}
where ${\rm {\bf b}}$ is the specific body force and $\left( {{\rm {\bf
\upsilon }}\cdot { }^\ast \nabla ^{\left( \beta \right)}} \right){\rm {\bf
\upsilon }}$ is the convection term of the fractional flow

From (\ref{eq22}) and (\ref{eq23}) the systems of the fractional Navier--Stoke equations are
given by
\begin{equation}
\label{eq24}
\left\{ {{\begin{array}{*{20}c}
 {\frac{\partial \rho }{\partial t}+{\rm {\bf \upsilon }}\cdot \left( {{
}^\ast \nabla ^{\left( \beta \right)}\rho } \right)=0,} \hfill \\
 {\rho \frac{\partial {\rm {\bf \upsilon }}}{\partial t}={ }^\ast \nabla
^{\left( \beta \right)}\cdot \left[ {-p{\rm {\bf I}}+2\mu {\rm {\bf \Lambda
}}+\lambda \left( {{ }^\ast \nabla ^{\left( \beta \right)}\cdot {\rm {\bf
\upsilon }}} \right){\rm {\bf I}}} \right]+\rho {\rm {\bf b}}-\rho \left(
{{\rm {\bf \upsilon }}\cdot { }^\ast \nabla ^{\left( \beta \right)}}
\right){\rm {\bf \upsilon }}£} \hfill \\
 {{\rm {\bf \upsilon }}={\rm {\bf \upsilon }}_0 ,} \hfill \\
\end{array} }} \right.
\end{equation}
which leads to
\begin{equation}
\label{eq25}
\left\{ {{\begin{array}{*{20}c}
 {\frac{\partial \rho }{\partial t}+{\rm {\bf \upsilon }}\cdot \left( {{
}^\ast \nabla ^{\left( \beta \right)}\rho } \right)=0,} \hfill \\
 {\rho \frac{\partial {\rm {\bf \upsilon }}}{\partial t}=-{ }^\ast \nabla
^{\left( \beta \right)}\cdot p{\rm {\bf I}}\mbox{+}{ }^\ast \nabla ^{\left(
\beta \right)}\cdot 2\mu {\rm {\bf \Lambda }}\mbox{+}{ }^\ast \nabla
^{\left( \beta \right)}\cdot \lambda \left( {{ }^\ast \nabla ^{\left( \beta
\right)}\cdot {\rm {\bf \upsilon }}} \right){\rm {\bf I}}+\rho {\rm {\bf
b}}-\rho \left( {{\rm {\bf \upsilon }}\cdot { }^\ast \nabla ^{\left( \beta
\right)}} \right){\rm {\bf \upsilon }},} \hfill \\
 {{\rm {\bf \upsilon }}={\rm {\bf \upsilon }}_0 .} \hfill \\
\end{array} }} \right.
\end{equation}
The constitutive equation of the incompressible Navier--Stokes fluid can be
written in the form:

\[{\rm {\bf J}}=-p{\rm {\bf I}}+2\mu {\rm {\bf \Lambda }},\mbox{ }{ }^\ast
\nabla ^{\left( \beta \right)}\cdot {\rm {\bf \upsilon }}=0,\]
since
\begin{equation}
\label{eq26}
\frac{\partial \rho }{\partial t}+{\rm {\bf \upsilon }}\cdot \left( {{
}^\ast \nabla ^{\left( \beta \right)}\rho } \right)=\frac{\partial \rho
}{\partial t}+{ }^\ast \nabla ^{\left( \beta \right)}\cdot \left( {{\rm {\bf
\upsilon }}\rho } \right).
\end{equation}
In this case, Eq.(\ref{eq25}) can be rewritten in the form:
\begin{equation}
\label{eq27}
\left\{ {{\begin{array}{*{20}c}
 {{ }^\ast \nabla ^{\left( \beta \right)}\cdot {\rm {\bf \upsilon }}=0,}
\hfill \\
 {\rho \frac{\partial {\rm {\bf \upsilon }}}{\partial t}=-{ }^\ast \nabla
^{\left( \beta \right)}\cdot p{\rm {\bf I}}\mbox{+}2\mu { }^\ast \nabla
^{\left( \beta \right)}\cdot {\rm {\bf \Lambda }}+\rho {\rm {\bf b}}-\rho
\left( {{\rm {\bf \upsilon }}\cdot { }^\ast \nabla ^{\left( \beta \right)}}
\right){\rm {\bf \upsilon }},} \hfill \\
 {{\rm {\bf \upsilon }}={\rm {\bf \upsilon }}_0 .} \hfill \\
\end{array} }} \right.
\end{equation}
or
\begin{equation}
\label{eq28}
\left\{ {{\begin{array}{*{20}c}
 {{ }^\ast \nabla ^{\left( \beta \right)}\cdot {\rm {\bf \upsilon }}=0,}
\hfill \\
 {\rho \frac{\partial {\rm {\bf \upsilon }}}{\partial t}=-{ }^\ast \nabla
^{\left( \beta \right)}\cdot p{\rm {\bf I}}\mbox{+}\mu { }^\ast \nabla
^{\left( {2\beta } \right)}{\rm {\bf \upsilon }}+\rho {\rm {\bf b}}-\rho
\left( {{\rm {\bf \upsilon }}\cdot { }^\ast \nabla ^{\left( \beta \right)}}
\right){\rm {\bf \upsilon }},} \hfill \\
 {{\rm {\bf \upsilon }}={\rm {\bf \upsilon }}_0 .} \hfill \\
\end{array} }} \right.
\end{equation}
In the case, Eq.(\ref{eq28}) is called as the fractional Navier--Stoke equations.

Without the specific body force in Eq.(\ref{eq28}), we obtain a new form of the
fractional Navier--Stoke equations
\begin{equation}
\label{eq29}
\left\{ {{\begin{array}{*{20}c}
 {{ }^\ast \nabla ^{\left( \beta \right)}\cdot {\rm {\bf \upsilon }}=0,}
\hfill \\
 {\rho \frac{\partial {\rm {\bf \upsilon }}}{\partial t}=-{ }^\ast \nabla
^{\left( \beta \right)}\cdot p{\rm {\bf I}}\mbox{+}\mu { }^\ast \nabla
^{\left( {2\beta } \right)}{\rm {\bf \upsilon }}-\rho \left( {{\rm {\bf
\upsilon }}\cdot { }^\ast \nabla ^{\left( \beta \right)}} \right){\rm {\bf
\upsilon }},} \hfill \\
 {{\rm {\bf \upsilon }}={\rm {\bf \upsilon }}_0 .} \hfill \\
\end{array} }} \right.
\end{equation}
\textbf{Conclusions }

In our work, we used the fractional gradient and Laplacian operators via
fractional derivatives without singular kernel to investigate the
mathematical theory for the heat and fluid flows. The fractional
heat-conduction equations and the fractional Navier-Stokes equations were
discussed. The results have opened the new directions of the heat and fluid
flows within fractional derivatives without singular kernel.

\textbf{References}

\begin{enumerate}
\item K. B Oldham J Spanier The fractional calculus: theory and applications of differentiation and integration to arbitrary order Academic Press, New York, 1974.
\item J. Sabatier, O. P. Agrawal, J. T. Machado, Advances in fractional calculus Springer, 2007.
\item A. A. Kilbas, H. M. Srivastava, J. J. Trujillo, Theory and Applications of Fractional Differential Equations, Academic Press, New York, 2006.
\item R. Gorenflo, F. Mainardi, Fractional calculus and stable probability distributions, Archives of Mechanics 50 (1998) 377-388.
\item V. E. Tarasov, Heat transfer in fractal materials, International Journal of Heat and Mass Transfer 93 (2016) 427-430.
\item Y. Z. Povstenko, Thermoelasticity that uses fractional heat conduction equation, Journal of Mathematical Sciences 162 (2009), 296-305.
\item M. A. Ezzat,. Thermoelectric MHD non-Newtonian fluid with fractional derivative heat transfer, Physica B 405 (2010) 4188-4194.
\item M. Khan, T. Hayat, S. Asghar, Exact solution for MHD flow of a generalized Oldroyd-B fluid with modified Darcy's law, International Journal of Engineering Science 44 (2006) 333-339.
\item M. Caputo, M. Fabrizio, A new definition of fractional derivative without singular Kernel,~Progress in Fractional Differentiation and Applications~1 (2015) 73-85.
\item Lozada, J.; Nieto, J. J. Properties of a new fractional derivative without singular kernel,~Progress in Fractional Differentiation and Applications 1 (2015) 87-92.
\item A, Alsaedi J. J Nieto, V. Venktesh, Fractional electrical circuits, Advances in Mechanical Engineering 7 (2015) 1-7.
\item X. J. Yang, H. M. Srivastava, J. A. Machado, A new fractional derivative without singular kernel: Application to the modelling of the steady heat flow, Thermal Science, DOI: 10.2298/TSCI151224222Y, arXiv:1601.01623, 2015.
\item A M. Yang Y. Han, J. Li, W. X. Liu, On steady heat flow problem involving Yang-Srivastava-Machado fractional derivative without singular kernel, Thermal Science 20 (2015) in press.
\item X. J. Yang, Advanced local fractional calculus and its applications, World Science, New York, 2012.
\item C. Cattani, H. M. Srivastava X. J. Yang, (Eds.) Fractional dynamics, Walter de Gruyter, 2016.
\item X. J. Yang, D. Baleanu, H. M. Srivastava, Local Fractional Integral Transforms and Their Applications, Academic Press, New York, 2015.
\end{enumerate}

\end{document}